\documentstyle[12pt,leqno]{article}
 \newcommand{\beq}{\begin{equation}}
\newcommand{\eeq}{\end{equation}}
 \newcommand{\beth}{\begin{theo}}
\newcommand{\eth}{\end{theo}}

\makeatletter%%%%%%%%%%%%
\@addtoreset{equation}{section}
\makeatother%%%%%%%%%%%%

\newcommand{\Subset}{\subset\subset}
\newcommand{\Psux}{{ \Psi_{u,x}}}
\newcommand{\Psvx}{{ \Psi_{v,x}}}
\newcommand{\Pspx}{{ \Psi_{\vph,x}}}
\newcommand{\Pspo}{{ \Psi_{\vph,0}}}
\newcommand{\Psuo}{{ \Psi_{u,0}}}
\newcommand{\nux}{\nu(u,x)}
\newcommand{\nuxa}{\nu(u,x,a)}
\newcommand{\nuph}{\nu(u,\vph)}
\newcommand{\mrph}{\mu_r^\vph}
\newcommand{\mrP}{\mu_r^\Phi}
\newcommand{\arP}{\gamma_r^\Phi}
\newcommand{\aoP}{\gamma_{-1}^\Phi}
\newcommand{\tmux}{T_{m,x}u}
\newcommand{\obl}{\Omega}
\newcommand{\nol}{\{0\}}
\newcommand{\okn}{1\le k\le n}
\newcommand{\vph}{\varphi}

\newcommand{\Rn}{{ \bf R}^n}
\newcommand{\Rnm}{{ \bf R}_-^n}
\newcommand{\Rnp}{{ \bf R}_+^n}
\newcommand{\Cn}{{ \bf C}^n}

\newtheorem{theo}{Theorem}
\newtheorem{cor}{Corollary}
\newtheorem{lem}{Lemma}
\newtheorem{prop}{Proposition}

\begin{document}

 \title{Lelong numbers with respect to regular plurisubharmonic weights}
\author{ Alexander RASHKOVSKII  }

\date{}
\maketitle

\begin{abstract}
Generalized Lelong numbers $\nu(T,\vph)$ due to Demailly are specified 
for the case of positive closed currents
$T=dd^cu$ and plurisubharmonic weights $\vph$ with multicircled
asymptotics. Explicit formulas for these values are obtained in terms
of the directional Lelong numbers of the functions
$u$ and the Newton diagrams of
$\vph$. An extension of Demailly's approximation theorem 
is proved as well.
\end{abstract}

\section{Introduction}

A standard quantative characteristic for singularity of 
a plurisubharmonic function $u$ at
a point $x\in\Cn$ is its Lelong number
$$
\nux=\lim_{r\to 0}\int_{|z-x|<r} dd^c u\wedge (dd^c\log|z-x|)^{n-1};
$$
here $d=\partial + \bar\partial,\ d^c= ( \partial
-\bar\partial)/2\pi i$. When $u=\log|f|$, $f$ being a holomorphic
function with $f(x)=0$, $\nux$ is just the multiplicity of the zero of 
$f$ at the point $x$. The Lelong number can also be calculated as
\beq
\label{eq:form1}
\nux=\lim_{r\to -\infty}r^{-1}\sup\{u(z):\: |z-x|\le e^r\}
= \lim_{r\to -\infty}r^{-1}{\cal M}(u,x,r),
\eeq
where ${\cal M}(u,x,r)$ is the mean value of $u$ over the sphere
$|z-x|=e^r$, see \cite{keyKis1}. Various results on Lelong numbers
and their applications to complex analysis can be found in
\cite{keyLe1}, \cite{keyLeG}, \cite{keyHo}, \cite{keyLe2}.

A more detailed information on the
behaviour of $u$ near $x$ can be obtained by means of the {\it refined},
or {\it directional, Lelong numbers} \cite{keyKis2}
\begin{eqnarray}
\nuxa &=& \lim_{r\to -\infty}r^{-1}\sup\{u(z):\: |z_k-x_k|\le
           e^{ra_k},\ \okn\}\nonumber\\
    &=& \lim_{r\to -\infty}r^{-1}\lambda(u,x,ra),
\label{eq:dir}
\end{eqnarray}
where $a=(a_1,\ldots,a_n)\in\Rnp$ and $\lambda(u,x,b)$ is the mean value
of $u$ over the set $\{z:\: |z_k-x_k|=\exp{b_k},\ \okn\}$.

A general notion of the Lelong number with respect to a plurisubharmonic weight
was introduced by J.-P.~Demailly \cite{keyD1}. Let $\vph$ be a
semiexhaustive plurisubharmonic function on a domain $\obl\subset\Cn$, that is, 
$B_R^\vph:=\{z:\: \vph(z)<R\}\Subset\obl$ for some real $R$. The value
\beq
\label{eq:gen}
\nuph=\lim_{r\to -\infty}\int_{B_r^\vph} dd^c u\wedge\left[ dd^c\vph
\right]^{n-1}
\eeq
is called {\it the generalized Lelong number} of $u$ with respect to
the weight $\vph$. For a detailed study of this notion, see
\cite{keyD3}. An analog of formula (\ref{eq:form1}),
if $ \vph$ satisfies $(dd^c\vph)^n=0$ on $\obl\setminus
\vph^{-1}(-\infty)$, is the relation 
$$
\nuph=\lim_{r\to -\infty}\mrph(u)
$$
where $\mrph$ is the swept out Monge-Amp\`ere measure for
$(dd^c\vph)^n$ on the pseudosphere $S_r^\vph:=\{z:\: \vph(z)=r\}$, i.e.
\beq
\mrph(u)=\int_{S_r^\vph} u\left[(dd^c\vph_r)^n-(dd^c\vph)^n\right]
\label{eq:swept}
\eeq
($\vph_r=\max\{\vph,r\}$). In particular,
$\nux=\nu(u,\log|\cdot-x|)$ and 
\beq
\label{eq:nuxa}
\nuxa=a_1\ldots a_n\,\nu(u,\vph_{a,x})
\eeq
with the weight
\beq
\vph_{a,x}(z)=\max_k \,a_k^{-1}\log|z_k-x_k|.
\label{eq:phax}
\eeq

Actually, the Lelong number $\nuph$ is a function of asymptotic behaviour
of the functions $u$ and $\vph$ near $x$, that follows from the Comparison
Theorem due to Demailly, which we state here in a  form
convenient for our purposes.

\medskip
{\bf Theorem A} (\cite{keyD1}, Th. 5.9).
{\sl Let $u_1$ and $u_2$ be plurisubharmonic functions 
 on a neighbourhood of a point $x\in\Cn$, $\vph_1$ and $\vph_2$
be plurisubharmonic weights with $\vph_1^{-1}(-\infty)=
\vph_2^{-1}(-\infty)=x$.
Suppose that $u_1(x)=-\infty$,
$$
\limsup_{z\to x}\frac{u_2(z)}{u_1(z)}\le 1
$$
and
$$
\limsup_{z\to x}\frac{\vph_2(z)}{\vph_1(z)}\le 1.
$$
Then}
$\nu(u_2,\vph_2)\le \nu(u_1,\vph_1)$.

\medskip
The generalized Lelong numbers give a powerful 
and supple instrument for investigation of singularities of 
plurisubharmonic functions.
Another thing is that one pays for the universality of 
such numbers
with lack of explicit ways for their evaluation. The objectives for the
present note are to look for a subclass of the weights which is wide
enough and at the same time convenient for treatment. Theorem~A suggests
that one should try to consider weights with certain "regular" asymptotics.
The condition $\vph(z)\sim\log|z-x|$ reduces the situation to the
standard Lelong numbers and gives nothing new. To deal with more
refined asymptotics,  we consider here two classes of the weights.

The first one uses the notion of local indicator \cite{keyLeR}.
Let a plurisubharmonic function $\Phi$ defined in the unit polydisk 
$D=\{z\in\Cn: |z_k|<1,\ 
\okn\}$
be nonpositive there and satisfy the relation
\beq
\label{eq:hom}
\Phi(z)=\Phi(|z_1|, \ldots,|z_n|)=
c^{-1}\Phi(|z_1|^c, \ldots,|z_n|^c)\quad\forall c>0.
\eeq
We will call such functions (abstract) {\it indicators}. The homogeneity
of indicators implies $(dd^c\Phi)^n=0$ outside the origin, provided
$\Phi^{-1}(-\infty)=0$. Such functions seem to be good candidates for
the weights we are looking for. Besides, the indicators present a scale of
plurisubharmonic characteristics for local behavior of plurisubharmonic
functions near their singularity points. Namely,
given a plurisubharmonic function $v$, its {\it local indicator} 
at a point $x$ is 
a plurisubharmonic function $\Psvx$ in the unit polydisk 
$D$ such that
\beq
\Psvx(y)=-\nu(v,x,a),\quad a=-(\log|y_1|,\ldots,\log|y_n|).
\label{eq:lind}
\eeq
It is the
largest negative plurisubharmonic function in $D$ whose directional
Lelong numbers
at $0$ coincide with those of $v$ at $x$, so 
\beq
v(z)\le\Psvx(z-x)+C
\label{eq:bound}
\eeq
near $x$. Besides, as was shown in \cite{keyR1}, $\Psvx$ can be
described as the limit (in $L_{loc}^1$) of the sequence 
\beq
\left(T_{m,x}v\right)(y)=m^{-1}v(x+y^m)
\label{eq:tg}
\eeq
as $m\to\infty$; here $y^m=(y_1^m,\ldots,y_n^m)$, $m\in{\bf Z}_+$.
Moreover, for a multicircled function $v$ negative in the unit polydisk,
the functions $v_R(z):=R^{-1}v(|z_1|^R,\ldots,|z_n|^R)$ increase to
a function $V(z)$ as $R\to+\infty$, and $V^*=\Psi_{v,0}$.

Local indicators are obviously indicators, and due to relations (\ref{eq:hom})
%Due to the obvious relations 
%\beq
%\Psvx(y)=\Psvx(|y_1|,\ldots,|y_n|)=c^{-1}\Psvx(|y_1|^c,\ldots,|y_n|^c)
%\quad\forall c>0,
%\label{eq:hom}
%\eeq
their use is quite efficient. 
Note that, in view of the definition of the local indicator and relation
(\ref{eq:nuxa}),
we have
\beq
\label{eq:i1}
\nu(u, \vph_{a,x})=\nu(\Psux,\vph_{a,x})=\nu(\Psux,\Psi_{\vph_{a,x}}).
\eeq
The second equation is evident since $\Psi_{\vph_{a,x}}={\vph_{a,x}}$,
while the first one is of some interest because it does not follow
from Theorem~A (no assumption on asymptotic behaviour of $u$ is made).

We will say that a weight $\vph$ with $\vph^{-1}(-\infty)=x$ is 
{\it almost homogeneous} if it is asymptotically equivalent to its indicator
$\Psi_{\vph,x}$, that is,
\beq
\label{eq:reg}
\exists\lim_{z\to x}{\vph(z)\over \Psi_{\vph,x}(z-x)}=1.
\eeq
It is easy to see that if the limit exists, it necessarily equals $1$.
Besides, the residual Monge-Amp\`ere measure of every 
almost homogeneous weight $\vph$
at $x$, $(dd^c\vph)^n(x)$, coincides with that of its indicator.
An example of such a weight is $\vph(z)=\log|F(z)|$ with a holomorphic
mapping $F:\obl\to{\bf C}^m$, $m\ge n$, $F^{-1}(0)=x$, such that
$$
\lim_{z\to x}{\log|F(z)| \over \sup\limits_{J\in\omega_x}\log|(z-x)^J|}=1,
$$
$\omega_x$ being the collection of all multi-indices $J$ satisfying
$\partial^J F/\partial z^J (x)\neq 0$. Generalized pluri-complex 
Green functions with respect to given indicators \cite{keyLeR} give
another example of almost homogeneous weights.

Note that every plurisubharmonic weight $\vph$ is the limit of a decreasing
sequence of almost homogeneous weights with the same local indicator
as $\vph$.
Indeed, these are $\vph_N(z)=\sup\{\vph(z),\Psi_{\vph,x}(z-x)-N\}$,
$N>0$.

Below we show that the generalized Lelong numbers with respect to 
almost homogeneous
weights inherit some nice properties from the standard and directional 
Lelong numbers. In particular, calculation of $\nuph$ can be reduced to 
that for 
the indicators, both of $u$ and $\vph$ (note that no regularity of $u$ 
is assumed). Namely,
let $\vph_x(z):=\vph(z-x)$ and $\vph^{-1}(-\infty)=\nol$,
then for any plurisubharmonic function $u$ in a domain $\obl\subset\Cn$,
\beq
\label{eq:res1}
\nu(u,\vph_x)=\nu(\Psux,\vph)=\nu(\Psux,\Pspo)\quad
\forall x\in\obl
\eeq
(Theorem \ref{theo:1}),
which is an extension of relation (\ref{eq:i1}). 
Further, by a slight modification of
Demailly's arguments \cite{keyD2} we show that any plurisubharmonic 
function $u$
in a bounded pseudoconvex domain $\obl\subset\Cn$ can be approximated
by a sequence of functions 
$$
u_m={1\over 2m}\log\sum_l |\sigma_{ml}|^2,
$$
$\sigma_{ml}$ being holomorphic in $\obl$, such that for every  
almost homogeneous weight $\vph$,
\beq
\label{eq:res2}
\nu(u_m,\vph_x)\le\nu(u,\vph_x)\le\nu(u_m,\vph_x)+{A\over m}
\quad\forall x\in\obl
\eeq
with some constant $A=A(\vph)$ (Theorem \ref{theo:3}).

Finally, we give a geometric description for the swept out
Monge-Amp\`ere measures $\mu_r^\Phi$ for indicators $\Phi$ 
(Theorem~\ref{theo:4}), 
which leads to explicit
formulas for the numbers $\nuph$ with almost homogeneous weights $\vph$
in terms of the directional Lelong numbers of $u$ and $\vph$ 
(Corollary~\ref{cor:1}). When $\vph=\log|g|$ with a holomorphic mapping
$g$, $g(0)=0$, this reduces to computation on the Newton diagram of $g$
at the origin. 

%%%%%%%
\medskip

Another choice of a class of the weights are
those  whose behaviour near $x$ is asymptotically 
independent of the arguments
of $z_k-x_k$, $\okn$. Namely, we will say that 
a plurisubharmonic function
$\vph$ on a domain $\Omega\subset\Cn$ has a {\it multicircled singularity}
at a point $x\in\Omega$ (or that $\vph$ is
{\it almost multicircled} near a point $x\in\Omega$)  
if there exists a multicircled plurisubharmonic function $\lambda$ 
 (i.e. $\lambda(z)=\lambda(|z_1|,\ldots,|z_n|)$) in a neighbourhood 
of the origin,
such that
\beq
\exists\lim_{z\to x}{\vph(z)\over \lambda(z-x)}=1.
\label{eq:amc}
\eeq
In the terminology of \cite{keyZa}, it means that $\vph$ has 
a standard singularity generated by a multicircled function.
It is easy to see that $\vph$ has multicircled singularity at $x$ 
if and only if it satisfies
relation (\ref{eq:amc}) with $\lambda$ equal to some "circularization" 
of $\vph$,
say, to the mean value
\beq
\lambda=\lambda_{\vph,x}(z) =(2\pi)^{-n}
\int_{[0,2\pi]^n}\vph(x_1+z_1e^{i\theta_1},
\ldots,
x_n+z_ne^{i\theta_n})\,d\theta
\label{eq:circ}
\eeq
or to its maximum on the same set. 
 
Evidently, every almost homogeneous weight is almost multicircled, 
however the converse
is not true. On the other hand, it can be seen that each multicircled
weight $\vph$ has the same residual Monge-Amp\`ere measure at $x$ as its
indicator has at the origin
\cite{keyR3}, so one might hope that the above results could be
extended to the whole class of almost multicircled weights. This turns to be 
really the case, but only when it concerns plurisubharmonic functions $u$
whose $-\infty$ sets do not contain lines parallel to the coordinate axes.
Namely, we prove (\ref{eq:res1}) 
and ({\ref{eq:res2}) 
for every almost
multicircled weight $\vph$ and functions $u$ with the above extra condition.
And, surprisingly, it fails to be true, for instance,
when $u(z_1,0,\ldots,0)\equiv -\infty$ and $x=0$.

\bigskip
{\it Notation.} Throughout the paper, $D$ is the unit polydisk in $\Cn$,
$\obl$ is a domain in $\Cn$, and 
$PSH(\obl)$ is the collection of all plurisubharmomic functions on $\obl$.
If $x\in\obl$ and  a function
$u\in PSH(\Omega)$
is such that its restriction to each line $\{z:z_j=x_j\ 
\forall j\neq k\}$, $\okn$, 
is not identically $ -\infty$, then we will write $u\in PSH_*(\Omega,x)$. 
The Lelong number of $u$ at $x$ is denoted $\nux$, and 
$\nuxa$ and $\nuph$ are its directional
(\ref{eq:dir}) and generalized Lelong numbers (\ref{eq:gen}).
Any plurisubharmonic function $\Phi$ in $D$ satisfying
%with $\Phi(z)=\Phi(|z_1|, \ldots,|z_n|)=
%c^{-1}\Phi(|z_1|^c, \ldots,|z_n|^c)$ for all $c>0$ 
(\ref{eq:hom}) will be called an indicator,
and the function $\Psvx$ defined by (\ref{eq:lind}) is the (local) indicator of
$v$ at $x$. 
The class of all plurisubharmonic weights $\vph$ in a neighborhood
of the origin, $\vph^{-1}(-\infty)=\nol$,
will be denoted by $W$, and its subclasses consisting of
almost homogeneous and almost multicircled weights,
in the sense of (\ref{eq:reg}) and (\ref{eq:amc}), will be denoted by $W^h$
and $W^m$, respectively. 

%%%%%%%%%%%%%%%%%%%%%%%%%%%%%%%

\section{Reduction to indicators}

%Let $\vph$ be an arbitrary plurisubharmonic weight in a neighborhood
%of the origin, $\vph^{-1}(-\infty)=\nol$. The class of all such
%weights will be denoted by $W$, and its subclasses consisting of
%regular and almost multicircled weights,
%in the sense of (\ref{eq:reg}) and (\ref{eq:amc}), will be denoted by $W^h$
%and $W^m$, respectively. Besides, let a plurisubharmonic function
%$u$ on a domain $\Omega$
%be such that its restriction to each line $\{z:z_j=x_j\ 
%\forall j\neq k\}$, $\okn$, 
%is not identically $ -\infty$, $x$ being a fixed point of $\Omega$. 
%The collection of all such functions
%$u$ will be denoted by $PSH_*(\Omega,x)$. 
%Note that if $u\in PSH_*(\Omega,x)$ then its indicator at $x$
%is locally bounded outside the origin.

\medskip
\begin{prop}
Let $u\in PSH(\obl)$, $x\in\obl$, and the functions $\tmux$ be defined 
by (\ref{eq:tg}). Then for each weight $\vph\in W$
 there exists the limit
$$
\lim_{m\to\infty} \nu(\tmux,\vph)=
\nu(\Psux,\vph).
$$
\label{prop:1}
\end{prop}

{\it Proof.} As was mentioned in Introduction, 
$\tmux\to\Psux$ in $L_{loc}^1(D)$ (\cite{keyR1}, Theorem~8). 
Semicontinuity theorem for
generalized Lelong numbers (Proposition~3.12 of \cite{keyD3}) then implies
$$
\limsup_{m\to\infty}\nu(\tmux,\vph)\le\nu(\Psux,\vph).
$$
On the other hand, the functions $\tmux$ have the same indicator
at the origin for all $m$, and it coincides with $\Psux$.
By (\ref{eq:bound}),
$\tmux\le\Psux +C_m$ near the origin, so 
$\nu(\tmux,\vph)\ge\nu(\Psux,\vph)$ $\forall m$, and the proof 
is complete.

\medskip

Now we specify the weight $\vph$ to be almost multicircled. For a function
$f$ defined on a subset of $\Cn$, $f_x(z)$ will denote $f(z-x)$,
$x\in\Cn$.  Then
\beq
\nu(u,\vph_x)=\nu(u_{-x},\vph)=\nu(u_{-x},\lambda)
\label{eq:dem}
\eeq
with multicircled $\lambda$ from (\ref{eq:amc}), 
the second equation being a consequence of relation (\ref{eq:amc})
in view of Theorem~A.
Actually, a stronger relation takes place.

\medskip
\beth
a) If $\vph\in W^m$, $x\in\Omega$, then for every function
$u\in PSH_*(\obl,x)$
\beq
\label{eq:1.1}
\nu(u,\vph_x)=\nu(\Psux,\vph)=\nu(\Psux,\Pspo);
\eeq

b) if $\vph\in W^h$ then (\ref{eq:1.1}) holds for every function
$u\in PSH(\obl)$;

c) there exist a multicircled function $u$ and a weight $\vph\in W^m$
such that $\nu(u,\vph)>\nu(\Psuo,\Pspo)$.

\label{theo:1}
\eth

{\it Proof}. When $\vph$ is an indicator (i.e., $\Psi_{\vph,0}=\vph$), 
relation (\ref{eq:1.1}) follows from Proposition~\ref{prop:1}:
\begin{eqnarray*}
\nu(\tmux,\vph) &= &\lim_{r\to -\infty}\int_{B_r^\vph}
                  dd^c\tmux\wedge (dd^c\vph)^{n-1} \\
 &=&\lim_{r\to -\infty}\int_{B_{mr}^\vph}
                  dd^c u_{-x}\wedge (dd^c\vph)^{n-1}
\end{eqnarray*}
by the homogeneity of indicators (see (\ref{eq:hom})).
The right-hand side equals $\nu(u_{-x},\vph)$, and the statement
follows from (\ref{eq:dem}) and Proposition \ref{prop:1}. 

By Theorem~A, this implies b).
 
To prove a), we need the following

\begin{lem}
\label{lem:1}
Let $v\in PSH_*(D,0)$ be multicircled in the unit polydisk, 
then for every $r\in (0,1)$ 
there exists a constant $A>0$ 
such that 
$$
v(z)\ge A\sup_j \log|z_j|\quad
$$
for all $z$, $|z_k|\le r$, $\okn$.
\end{lem}

{\it Proof of Lemma \ref{lem:1}.} 
Consider the function $v_1(\zeta)=v(\zeta,0,\ldots,0)-C_1$
with $C_1$ such that $\sup \{v_1(\zeta):|\zeta|<1\}=0$.
Since the ratio $v_1(\zeta)/\log|\zeta|$ decreases to $\nu_1\ge 0$ as
$|\zeta|\searrow 0$, we have $v_1(\zeta)\ge A_1\log|\zeta|$ for some $A_1>0$
and all $\zeta$ with $|\zeta|\le r$. Therefore,
$v(z)\ge v_1(z_1) \ge A_1\log|z_1|+C_1\ge A'_1\log|z_1|$, $|z_1|<r$.
The same arguments for $j=2,\ldots,n$ complete the proof of the lemma.

\medskip
Let now $\vph\in W^m$ and $u\in PSH_*(\obl,0)$. 
By (\ref{eq:dem}) we may assume 
$\vph$ to be multicircled. In this case, $\nu(u,\vph_x)=\nu(\lambda_{u,x},\vph)$
with the function $\lambda_{u,x}$ defined by (\ref{eq:circ}). So, it suffices
to prove the assertion for nonpositive,
multicircled functions $u\in PSH_*(D_{2r},0)$
in a polydisk $D_{2r}=\{z:|z_k|<2r,\ \okn\}$, $r\in (0,1/2)$.

%%%%%%%%%%%%
As was mentioned in 
Introduction, the functions 
$$
u_R(z)=R^{-1}u(|z_1|^R,\ldots,|z_n|^R)
\nearrow U(z)
$$ 
and 
$$
\vph_R(z)=R^{-1}\vph(|z_1|^R,\ldots,|z_n|^R)
\nearrow \Phi(z)
$$
as $R\to +\infty$, with $U^*=\Psuo$ and $\Phi^*=\Pspo$.
Let $A>0$ be chosen as in Proposition~\ref{prop:1}
for both $u$ and $\vph$, 
and $L=2A\log r$.
Denote $v_R(z)=\max\{u_R(z),L\}$,
$w_R(z)=\max\{\vph_R(z),L\}$, $V(z)=\max\{U(z),L\}$,
 $W(z)=\max\{\Phi(z),L\}$, $P(z)=\max\{\Psuo(z),L\}$, 
$Q(z)=\max\{\Pspo(z),L\}$.
By the choice of $L$, $v_R=u_R$ and $w_R=\vph_R$ near $\partial D_r$ 
for all
$R\ge 1$, as well as $P=\Psuo$ and $Q=\Pspo$ there. Then
\beq
\label{eq:3a}
\int_{D_r} dd^c v_R\wedge
 (dd^cw_R)^{n-1}=\int_{D_r}dd^c u_R\wedge (dd^c
\vph_R)^{n-1}\ge\nu(u_R,\vph_R)=\nu(u,\vph)
\eeq
since
$$
\nu(u_R,\vph_R)=\lim_{r\to 0}\int_{D_r} dd^c u_R\wedge(dd^c\vph_R)^{n-1}=
\lim_{r\to 0}\int_{D_{r^R}}dd^cu\wedge (dd^c\vph)^{n-1}=\nu(u,\vph).
$$

On the other hand, $v_R\nearrow V$, $w_R\nearrow V$, $V^*=P$, and $W^*=Q$, 
so by the  
convergence theorem for increasing sequences of bounded
plurisubharmonic functions \cite{keyBT2},
$$
\lim_{R\to \infty}\int_{D_r} dd^cv_R\wedge(dd^cw_R)^{n-1}=
\int_{D_r} dd^cP\wedge(dd^cQ)^{n-1}=
\int_{D_r} dd^c\Psuo\wedge(dd^c\Psi_\vph)^{n-1}.
$$
Being compared with (\ref{eq:3a}) it gives us the relation
$\nu(u,\vph)\le\nu(\Psux,\Pspo)$.
As the opposite inequality is true due to Theorem~A, the proof 
of a) is complete.

Finally, 
consider the function $\vph(z_1,z_2)=\max\{-|\log|z_1||^{1/2},\log|z_2|\}$. 
Clearly, $\Pspo\equiv 0$. At the same time, for the function
$u(z)=\log|z_1|$, we have $\nu(u,\vph)=1$.

\medskip
{\it Remark.} If $u$ is a semiexhaustive plurisubharmonic
function in $\obl$, $u^{-1}(-\infty)=x$, its residual Monge-Amp\`ere measure 
$(dd^cu)^n(\{x\})$ at $x$ is just its generalized Lelong number $\nuph$
with respect to the weight $\vph=u$, while the residual measure of its
indicator $\Psux$ is, by the definition, the Newton number of $u$ at $x$  
\cite{keyR1}. 
Theorem~\ref{theo:1} then implies in particular
that for every almost multicircular function
$u\in PSH_*(\obl,x)$, its residual measure equals its Newton number, the
result proved earlier in \cite{keyR3}.

%%%%%%%%%%%%%%%%%%%%%%%%%%%%%%%%

\section{Approximation theorems}

An important theorem on approximation of plurisubharmonic functions was obtained
by J.-P.~Demailly (\cite{keyD2}, Proposition~3.1). 
Let $\Omega$ be a bounded pseudoconvex domain, $u\in PSH(\Omega)$,
and $\{\sigma_{ml}\}_l$ be an orthonormal basis of the Hilbert space
$$
H_m:=H_{m,u}(\Omega)=\{f\in Hol(\Omega):
\int_\Omega |f|^2e^{-2mu}\beta_n<\infty\}.
$$
Consider the functions
\beq
\label{eq:appr1}
u_m={1\over 2m}\log\sum_l |\sigma_{ml}|^2\in PSH(\Omega).
\eeq
Then there are constants $C_1, C_2>0$ such that for any point
$z\in\Omega$ and every $r<{\rm dist}\,(z,\partial \Omega)$,
$$
u(z)-{C_1\over m} \le u_m(z)\le \sup_{\zeta\in B_r(z)}u(\zeta)
+{1\over m}\log {C_2\over r^n}.
$$
In particular, $u_m\to u$ pointwise and in $L^1_{loc}(\Omega)$, and
$$
\nux-{n\over m}\le\nu(u_m,x)\le\nux\quad\forall x\in\Omega.
$$

Here we will show that actually the singularities of the functions
$u_m$ are almost the same that of $u$ not only in the sense of the
standard Lelong numbers but also with respect to arbitrary 
almost homogeneous weights,
as well as with almost multicircled weights (subject to the 
same restriction as in
Theorem~\ref{theo:1}).

We start with the following  modification of Demailly's result.

\beth
\label{theo:2}
Given a bounded pseudoconvex domain $\obl$,
there are constants $C_1, C_2>0$ such that for any function $u\in PSH(\obl)$
 and
every $z\in\obl$, the functions $u_m$ defined by (\ref{eq:appr1})
satisfy the relations
\beq
u(z)-{C_1\over m} \le u_m(z)\le \sup_{\zeta\in D_r(z)}u(\zeta)
+{1\over m}\log {C_2\over r_1\ldots r_n}
\label{eq:d1}
\eeq
for all $r=(r_1,\ldots, r_n)$ such that
$D_r(z)=\{\zeta: |z_k-\zeta_k|<r_k,\ \okn\}\Subset\obl$,
and 
\beq
\Psux(y)\le\Psi_{u_m,x}(y)\le\Psux(y)-m^{-1}\log|y_1\ldots y_n|
\quad\forall x\in\Omega,\ \forall y\in D.
\label{eq:d2}
\eeq
\eth

{\it Proof}. The first inequality in (\ref{eq:d1}) is the same as in 
 Demailly's Approximation theorem, and we repeat its proof here for 
completeness. Let $u(z)\neq -\infty$. By the Ohsawa-Takegoshi
extension theorem \cite{keyOT} applied to the $0$-dimensional subvariety
$\{z\}\subset\obl$, for any $a\in{\bf C}$ there exists a holomorphic
function $f$ on $\obl$ such that $f(z)=a$ and
$$
\int_\obl |f|^2\exp\{-2mu\}\,dV\le C|a|^2exp\{-2mu(z)\}
$$
with a constant $C=C(n, {\rm diam}\,\obl)$. Choosing $a$ such that
the right-hand side equals $1$ we get $f\in B_m$, the unit ball in the 
space $H_m$. Since
\beq
u_m(\zeta)=\sup_{g\in B_m}\,{\log|g(\zeta)|\over m}\quad
\forall\zeta\in\obl,
\label{eq:norm}
\eeq
it gives us 
$$
u_m(z)\ge{\log|f(z)|\over m}={\log|a|\over m}=u(z)-
{\log C\over 2m}.
$$

The proof of the second inequality in (\ref{eq:d1}) is a slight
modification of the corresponding arguments from \cite{keyD2}.
For any $g\in H_m$,
\begin{eqnarray*}
|g(z)|^2 &\le& {1\over\pi^n r_1^2\ldots r_n^2}
               \int_{D_r(z)} |g(\zeta)|^2\,dV\\
  &\le & {1\over\pi^n r_1^2\ldots r_n^2}\exp\{2m
\sup_{\zeta\in D_r(z)} u(\zeta)\}\,
               \int_{D_r(z)} |g|^2\exp\{-2mu\}\,dV,
\end{eqnarray*}
and so by (\ref{eq:norm})
$$
u_m(z)\le \sup_{\zeta\in D_r(z)} u(\zeta)+{1\over m}
 \log{\pi^{n/2}\over r_1\ldots r_n}.
$$

The first inequality in (\ref{eq:d1}) implies $\Psux\le\Psi_{u_m,x}$
$\forall x\in\obl$. To get the other bound, take any $y\in D$,
$y_1\ldots y_n\neq 0$, and $R>0$. Then for $r_k=|y_k|^R<
{\rm dist}\,(x,\partial\obl)$,
$$
{1\over R}\sup_{D_r(x)}u_m(\zeta)\le
{1\over R}\sup_{D_{2r}(x)}u_m(\zeta)-
{1\over m}{\log|y_1\ldots y_n|}+{C\over Rm},
$$
and the limit transition as $R\to\infty$ gives us the desired
inequality in view of the definition of the indicator
(\ref{eq:lind}) and (\ref{eq:dir}).

\medskip
{\it Remark}. In terms of the directional Lelong numbers, relations (\ref{eq:d2})
have the form
$$
\nuxa \le \nu(u_m,x,a)\le\nuxa
+m^{-1}\sum_j a_j
\quad\forall x\in\obl,\ \forall a\in \Rnp,
$$
so 
\beq
\label{eq:d10}
\lim_{m\to\infty}\nu(u_m,x,a)=\nuxa\quad \forall a\in \Rnp.
\eeq 
For the indicators,
the relation corresponding to (\ref{eq:d10}) 
is true for all $y$ with $y_1\ldots
y_n\neq 0$, while when $y_1\ldots y_n=0$ the regularization is needed:
$$
\Psux(y)=\limsup_{y'\to y}\,\lim_{m\to\infty}\Psi_{u_m,x}(y').
$$

\beth
\label{theo:3} In the conditions of Theorem \ref{theo:2},
$$
\nu(u_m,\vph_x)\le\nu(u,\vph_x)\le\nu(u_m,\vph_x)
+m^{-1}\sum_{\okn}\tau_k(\vph)\quad\forall x\in\obl
$$
for any almost homogeneous weight $\vph$ 
and $\tau_k(\vph)=\nu(\log|z_k|,\vph)$.

The same is true for every weight $\vph\in W^m$ provided $u\in PSH_*(\obl,x)$.
\eth

{\it Proof}. Let $\vph\in W^h$. 
Denote $\Phi=\Pspo$. By Theorems \ref{theo:1} and
\ref{theo:2},
$$
\nu(u_m,\vph_x)= \nu(\Psi_{u_m,x},\Phi)\le
\nu(\Psux,\Phi)=\nu(u,\vph_x).
$$
Similarly,
\begin{eqnarray*}
\nu(u,\vph_x) &=& \nu(\Psi_{u,x},\Phi)\le
       \nu(\Psi_{u_m,x},\Phi)+m^{-1}\sum_{\okn}\nu(\log|y_k|,\Phi)\\
&= &\nu(u_m,\vph_x)
+m^{-1}\sum_{\okn}\tau_k(\vph).
\end{eqnarray*}

If $u\in PSH_*(\obl,x)$, then (\ref{eq:d2}) implies $u_m\in PSH_*(\obl,x)$
for each $m$, so all the above arguments work
with arbitrary weights $\vph\in W^m$, too. The only exception is
the relation $\nu(\log|y_k|,\Phi)=\tau_k(\vph)$ which is to be replaced
with the inequality $\nu(\log|y_k|,\Phi)\le\tau_k(\vph)$.

%%%%%%%%%%%%%%%%%%%%%%%%%%%%%%%%%%%%%%%%%%%%%%%%%%%%%%%%%%%%%%%%%%%%%%%

\section{Swept out Monge-Amp\`ere measures}

The role of the swept out Monge-Amp\`ere measures $\mrph$ 
(\ref{eq:swept}) is demonstrated by the
Lelong-Jensen-Demailly formula \cite{keyD3}: 
if $u$ is a plurisubharmonic function in the pseudoball $B_R^\vph$, then for any
$r<R$,
$$
\mrph(u)-\int_{B_r^\vph} u(dd^c\vph)^n=
\int_{-\infty}^r\left[\int_{B_t^\vph}dd^c u\wedge (dd^c\vph)^{n-1}
\right]\,dt.
$$
When $(dd^c\vph)^n=\tau\,\delta_x$, $\delta_x$ being the Dirac
$\delta$-function at $x$, the Lelong-Jensen-Demailly formula
leads to the representation formula for plurisubharmonic functions:
$$
u(x)=\tau^{-1}\mrph(u)+\tau^{-1}\int_{-\infty}^r\left[\int_{B_t^\vph}
dd^c u\wedge (dd^c\vph)^{n-1}\right]\,dt,
$$
and $\nuph=\lim_{r\to -\infty} r^{-1}\mrph(u)$.

In the case of $\vph=\vph_{a,x}$ given by (\ref{eq:phax}),
$\mrph=(a_1\ldots a_n)^{-1}m_{ra}$ with $m_{ra}$ the normalized
Lebesgue measure on $\{|z_k|=\exp\{ra_k\},\ \okn\}$. No explicit
formulas are available in the general situation. However when studying
singularities, one interests mainly in asymptotic behavior of $\mrph$
as $r\to-\infty$. For regular weights $\vph$ it means that we may
restrict ourselves to study of the swept out measures for their
indicators $\Phi=\Pspx$.

Since $(dd^c\Phi)^n=0$ on $D\setminus\nol$,
$\mrP=(dd^c\Phi_r)^n$ for each $r<0$. The function $\Phi_r$ is
invariant under the rotations 
$$
(z_1,\ldots,z_n)\mapsto (z_1e^{i\omega_1},\ldots,z_ne^{i\omega_n}),
$$
and so is $\mrP$. Therefore we can write
$$
\mrP=(2\pi)^{-n}d\theta\otimes d\rho_r^\Phi
$$
with some measure $\rho_r^\Phi$ defined on the set
 $\{a\in\Rnp:\Phi(a)=r\}$.
Moreover, since $\mrP$ has no masses on the pluripolar set
$S_r^\Phi\cap\{z:z_1\ldots z_n=0\}$, we can pass to the coordinates
$z_k=\exp\{t_k+i\theta_k\}$, $\okn$. The functions
$$
f(t)=\Phi(e^{t_1},\ldots,e^{t_n})
$$ 
and
$$
f_r(t)=\Phi_r(e^{t_1},\ldots,e^{t_n})= \max\,\{f(t),r\}
$$
are convex in $\Rnm=-\Rnp$ and increasing in each $t_k$.
Simple calculations show that in these coordinates $\rho_r^\Phi$
transforms to 
$$
\arP=n!\,{\cal MA}[f_r]
$$ 
where ${\cal MA}$ is the
{\it real} Monge-Amp\`ere operator, see the details in
\cite{keyR1}. We recall that for smooth functions $v$,
$$
{\cal MA}[v]=\det\left({\partial^2 v\over\partial t_j\partial t_k}
\right)\,dt,
$$
and it can be extended as a positive measure to any convex function
(see \cite{keyRaT}).
So,
\begin{eqnarray*}
\mrP(u) &=& \int_{[0,1]^n}(2\pi)^{-n}\int_{[0,2\pi]^n} u(z_1e^{i\theta_1},
     \ldots, z_ne^{i\theta_n})\,d\theta\,d\rho_r^\Phi(|z_1|,\ldots,|z_n|)\\
&=&n!\,\int_{\Rnm}\lambda(u,0,t)\,d\arP(t)
\end{eqnarray*}
($\lambda(u,0,t)$ is the mean value of $u$ over $\{|z_k|=e^{t_k},\
\okn\}$). 

Since $f_r(t)=|r|f_{-1}(t/|r|)$,
$$
\mrP(u)=n!\,\int_{\Rnm}\lambda(u,0,|r|t)\,d\aoP(t).
$$
So, we only have to find an explicit expression for the measure
$\aoP$.

First of all, ${\rm supp}\, \aoP\subset L^\Phi$ with
\beq
\label{eq:LPhi}
 L^\Phi=\{t\in\Rnm:
f(t)=-1\}.
\eeq
If $\Phi=\Psux$, then $L^\Phi=-\{b\in\Rnp:\nu(u,x,b)=1\}$.

 As was shown in \cite{keyRaT}, for any convex function
$v$ in a domain $G\subset\Rn$,
\beq
\int_{F} {\cal MA}[v]={\rm Vol}\,\omega(F,v)\quad\forall F\subset G,
\label{eq:u3}
\eeq
where
$$
\omega(F,v)=\bigcup_{t^0\in F} \{a\in\Rn:\: v(t)\ge v(t^0)+
\langle a,t-t^0\rangle\quad\forall t\in G\}
$$
is the gradient image of the set $F$ for the surface $\{y=v(x),\
x\in G\}$.

Given a subset $F$ of $L^\Phi$, we put
\beq
\Gamma_F^\Phi=\{a\in\Rnp:\: \sup_{t\in F}\langle a,t\rangle =
\sup_{t\in L^\Phi}\langle a,t\rangle=-1\}
\label{eq:gf}
\eeq
and
\beq
\Theta_F^\Phi=\{\lambda a:\: 0\le\lambda\le 1,\ a\in \Gamma_F^\Phi\}.
\label{eq:gth}
\eeq
If $\Phi=\Psux$, then
$$
\Theta_{L^\Phi}^\Phi=\{a\in\Rnp:\sup_b[\nu(u,x,b)-\langle
a,b\rangle]\ge 0\}.
$$

Note that $\Gamma_{L^\Phi}^\Phi$ is an unbounded convex subset of 
$\Rnp$ and $f$ is the restriction of its supporting function to $\Rnm$.
When $\Phi$ is the indicator of $\log|g|$ for a holomorphic mapping
$g=(g_1,\ldots,g_n)$, $g(0)=0$, the set $\Gamma_{L^\Phi}^\Phi$ is the Newton
diagram for $g$ and $\Rnp\setminus\Theta_{L^\Phi}^\Phi$ is the
Newton polyhedron for $g$ at $0$ \cite{keyR1} as defined in \cite{keyKu1}. 
In this case, $\Gamma_F^\Phi$ is the union of bounded faces of the polyhedron 
corresponding to $F$.

\begin{prop}
\label{prop:2} 
For any compact subset $F$ of $L^\Phi$,
$\Theta_F^\Phi=\omega(F,f_{-1})$.
\end{prop}

{\it Proof}.
If $a\in\omega(F,f_{-1})$ then for some $t^0\in F$,
\beq
\langle a,t^0\rangle\ge\langle a,t\rangle-f_{-1}(t)-1\quad\forall
t\in\Rnm.
\label{eq:th1}
\eeq
In particular,
\beq
\langle a,t^0\rangle\ge\langle a,t\rangle\quad\forall t\in L^\Phi.
\label{eq:th2}
\eeq
When $t\to 0$, (\ref{eq:th1}) implies $\langle a,t^0\rangle\ge -1$.
In view of (\ref{eq:th2}) it means that $a\in\Theta_F^\Phi$.

Let now $a=\lambda a^0,\ a^0\in\Gamma_F^\Phi,\ 0\le\lambda\le 1$. 
Then there is a point $t^0\in F$ such that
$$
\langle a,t^0\rangle=\sup_{t\in F}\langle a,t\rangle =
\sup_{t\in L^\Phi}\langle a,t\rangle=-\lambda.
$$
For any $t\in\Rnm$, $t/|f(t)|\in L^\Phi$. If $f(t)\le -1$, then
$$
\langle a,t^0\rangle\ge \langle a,t/|f(t)|\rangle\ge\langle
a,t^0\rangle=\langle a,t\rangle-f_{-1}(t)-1.
$$
If $f(t)=-\delta>-1$, then $t/\delta\in L^\Phi$ and
\begin{eqnarray*}
\langle a,t\rangle-f_{-1}(t)-1&=&\delta\langle a,t/\delta\rangle-1
+\delta\le\delta\,\sup_{s\in L^\Phi}\langle a,s\rangle -1+\delta\\
&=&\delta\langle a,t^0\rangle-1+\delta =-\delta\lambda-1+
\delta\le\lambda=\langle a,t^0\rangle.
\end{eqnarray*}
The proposition is proved.

\begin{prop}
\label{prop:3}
The measure $\aoP$ is supported by $E^\Phi$, the set 
of extreme points of the convex set
$\{t: f(t)\le -1\}$.
\end{prop}

{\it Proof}. As
$$
\sup_{t\in L^\Phi}\langle a,t\rangle =
\sup_{t\in E^\Phi}\langle a,t\rangle\quad\forall a\in\Rnp,
$$
$\Theta_{L^\Phi}^\Phi=\Theta_{E^\Phi}^\Phi$. Hence
$\aoP(L^\Phi)=\aoP(E^\Phi)$ and thus $\aoP(F)=0$ for
every $F\Subset L^\Phi\setminus E^\Phi$.

\medskip

We have thus obtained the following

\beth
\label{theo:4}
For any plurisubharmonic function $u$ in a neighborhood of the origin
and for any indicator weight $\Phi$, the swept out Monge-Amp\`ere measure
$\mrP(u)$ on the set $\{\Phi(z)=r\}$, $r<0$, is determined by the formula
$$
\mrP(u)=n!\,\int_{E^\Phi}\lambda(u,0,|r|t)\,d\aoP(t)
$$
where $\lambda(u,0,|r|t)$ is the mean value of $u$ over the 
distinguished boundary of the polydisk $\{|z_k|<\exp\{|r|t_k\},\
\okn\}$ and
the measure $\aoP$ on the set $E^\Phi$ of extreme points of
the convex set
$\{t\in\Rnm: \Phi(e^{t_1},\ldots,
e^{t_n})\le-1\}$
is given by the relation $\aoP(F)={\rm Vol}\,\Theta_F^\Phi$ for
compact subsets $F$ of $E^\Phi$, $\Theta_F^\Phi$ being defined by
(\ref{eq:gf}), (\ref{eq:gth}), and (\ref{eq:LPhi}).
\eth

\begin{cor}
\label{cor:1}
If $\vph$ is an almost homogeneous weight, $\vph^{-1}(-\infty)=\{x\}$,
then for any plurisubharmonic function $u$ near $x$,
$$
\nuph=n!\,\int_{E^\Phi}\nu(u,x,-t)\,d\aoP(t)
$$
with $\Phi=\Pspx$ and the measure $\aoP$ the same as in 
Theorem~\ref{theo:4}.
If $u\in PSH_*(\obl,x)$ then the formula takes place
with arbitrary almost multicircled weight $\vph$.
\end{cor}

So, Corollary \ref{cor:1} gives a quantative expression for the fact that
the generalized Lelong number of a plurisubharmonic function $u$ with respect
to a regular weight $\vph$ is completely determined by the directional 
Lelong numbers of $u$ and $\vph$. 

Note that, under no regularity
condition on a weight $\vph$,
we have always the inequality
$$
\nuph\ge n!\,\int_{E^\Phi}\nu(u,x,-t)\,d\aoP(t).
$$

Of course, when $\vph=\log|g|$ with $g$ a holomorphic mapping, 
$g(0)=0$, the set $E^\Phi$ is finite and the measure $\aoP$
charges it with the volumes of the cones generated by the corresponding
$(n-1)$-dimensional faces of the Newton polyhedron of $g$.

%%%%%%%%%%%%%%%%%%%%%%%%%%%%%%%%%%
\bigskip

{\small{\it Acknowledgement.} The starting point of this work was
a question by Aydin Aytuna on explicit formulas for swept out 
Monge-Amp\`ere measures.
The author thanks Urban Cegrell and Vyacheslav Zakharjuta for
valuable discussions. The work was partially supported by the grant
INTAS-99-00089. 
 Part of the results were obtained during the author's
visit to the Feza Gursey Institute (Istanbul) under a NATO-T\"UBITAK grant.}

%%%%%%%%%%%%%%%%%%%%%%%%%%%%%%%%%%%%%%%%%%%%%%%%%%%%%%%%%%%%%%%%%%%%%%%%%%%%%%%%

\vskip 0.5cm

Mathematical Division, Institute for Low Temperature Physics
\par
47 Lenin Ave., Kharkov 61164,
Ukraine
%\par

\vskip0.1cm

E-mail: \quad rashkovs@ilt.kharkov.ua
% \quad rashkovs@ilt.kharkov.ua

%Phone: \quad +380 (572) 308-566

%Fax: \qquad \,+380 (572) 322-370

\end{document}